# Transforming Ridesharing: Harnessing Role Flexibility and HOV Integration for Enhanced Mobility Solutions*

Fatemeh Amerehi[1[0000-0002-6255-4573]] and Patrick Healy[1[0000-0002-3824-7442]]

[1] University of Limerick, Limerick, V94 T9PX, Ireland
Fatemeh.Amerehi@ul.ie

**Abstract.** While dynamic ridesharing has been extensively studied, there remains a significant research gap in exploring role flexibility within the many-to-many ridesharing scheme, where the system allows for several pickups for drivers and multiple transfers for riders. Previous works have predominantly assumed that all participants own a car and have focused on one-to-one arrangements. Additionally, there is a scarcity of research on integrating High Occupancy Vehicle (HOV) lanes and mathematical modelling. This study addresses these gaps by presenting a novel Mixed Integer Linear Programming (MILP) model that allows for role flexibility irrespective of car ownership and considers the implications of HOV lanes. Computational analysis highlights the benefits of incorporating role flexibility and accommodating non-car-owning participants in many-to-many ridesharing systems. Yet, excessive role shifts may create imbalances, impacting service to non-car owners. Further research should explore these correlations.

**Keywords:** Dynamic Ridesharing, Mixed Integer Linear Programming, Ride Matching, Role Flexibility, High Occupancy Vehicle

## 1 Introduction

Urban traffic congestion presents several challenges, including longer travel times, increased fuel consumption, higher emissions, and adverse environmental and health impacts. Expanding infrastructure and investing in public transit can be expensive and sometimes impractical in densely populated areas. While congestion pricing and alternative transportation modes exist, they may fall short in less densely populated areas. High-Occupancy Vehicle (HOV) lanes, dedicated for multi-occupant vehicles [1], effectively reduce congestion by promoting ridesharing and public transit. These lanes enhance traffic flow and serve as key strategies to tackle this issue, especially on busy freeways [1, 2]. HOV lanes in places like Beijing, Hong Kong, and California have shown significant benefits, including reduced travel times and improved traffic flow [3]. Carpoolers also enjoy incentives like reduced toll charges, motivating more people

* This work was conducted with the financial support of the Science Foundation Ireland Centre for Research Training in Artificial Intelligence under Grant No. 18/CRT/6223

to use high-occupancy transportation modes and thereby reducing overall road congestion [4].

Dynamic ridesharing systems with non-employee or non-contractor drivers also aim to alleviate congestion. They facilitate connections between travelers with compatible itineraries and schedules on short notice. Within ridesharing systems, an array of matching and transfer strategies exist. Some studies focus on one-to-one matching, where each rider is paired with a single driver capable of accommodating one passenger at a time [5]. A more comprehensive ridesharing approach involves a many-to-many scheme, where the system accommodates multiple stops for drivers and potential transfers between drivers and riders [1]. Another aspect of ridesharing systems is the role of participants as being driver or rider, which are typically defined upon trip registration. While some studies consider role flexibility [5, 6], the assumption is that all participants possess personal vehicles and thus have no preference, or that the system adheres to one-to-one frameworks. Similarly, studies on HOV [4] lack flexibility in participants roles, focusing on hired drivers and simplified one-to-one schemes. This study aims to fill these gaps by introducing a Mixed Integer Linear Programming (MILP) model for HOV lanes with a many-to-many approach and considering comprehensive role flexibility among participants. Ridesharing with multiple passengers in a vehicle may offer cost savings, and when combined with HOV lanes, it provides time-saving benefits.

## 2 Problem Description and Mathematical Formulations

We consider a dynamic ridesharing system that operates on a transportation network $\mathbb{G} = (\mathbb{N}, \mathbb{E})$, where $\mathbb{N}$ denotes nodes and $\mathbb{E}$ represents edges. Participants $\mathbb{P} = \{\mathbb{C} \cup \mathbb{NC}\}$ register their trip details and preferences are divided into car owners $\mathbb{C}$, and non-car owners $\mathbb{NC}$. Given the role-flexibility of the problem, individuals who take on the role of drivers $d$ form a subset $\mathbb{D} \subseteq \mathbb{C}$. Meanwhile, those who register as non-car owners or car owners but are assigned passenger roles constitute the set of $r$ riders $\mathbb{R} \supseteq \mathbb{NC}$. Each participant's travel data includes submission time $s_p$, origin $o_p$ and end $e_p$ locations, earliest departure $ED_p$, latest arrival times $LA_p$, with car owners specifying car capacity $Q_d$ and a notification deadline $\Delta_p$ for potential matches. We represent HOV nodes as $\mathbb{N}_H \subset \mathbb{N}$ and define $\mathbb{E}_H$ as the edges present in the sub-graph of $\mathbb{G}$ induced by $\mathbb{N}_H$.

To adapt to the dynamic nature of the problem and accommodate continuous information arrival we employ a rolling horizon strategy [5]. This involves periodic re-optimization of a static problem formulation at specific times based on incoming trip announcements. Prior to model solution, we perform preprocessing to account for users' spatiotemporal constraints. This narrows down the available nodes and links for each participant, leading to a reduced search space and facilitating the identification of viable and non-viable pairings of drivers and riders. We define the set of accessible nodes for participant $p \in \mathbb{P}$ as $\mathbb{AN}_p = \left\{i \in \mathbb{N}: T_{s_p,i} + T_{i,e_p} \leq LA_p - ED_p\right\}$, with $T_{i,j}$ denoting the shortest trip time between locations $i$ and $j$. Also, we establish sets of mutually accessible nodes for each participants-driver pair, denoted as $\mathbb{AN}_{pd}$, and their corresponding accessible links $\mathbb{AE}_{pd}$. Viable matches $\mathbb{RD}$ between riders and drivers, are based on overlapping accessible links and time constraints. It is important to note that due to transfer allowances a complete path or time overlap is not necessary. A match exists if $\mathbb{AE}_{pd} \neq \emptyset$ and $\neg(LA_d \leq ED_p \parallel LA_p \leq ED_d)$. To account for potential wait times, we introduce constraints into the model that synchronize meeting times at pickup and transfer points, all the while maintaining onboard riders' schedules. To formulate this problem, we define the decision variables as outlined in Table 1.

**Table 1.** Decision Variables

| Variable | Description |
|---|---|
| $x_{i,j}^{p,d}$ | Binary variable if driver $d$ serves participant $p$ from $i$ to $j$ |
| $y_i^{p,d}$ | Binary variable if driver $d$ picks up participant $p$ at node $i$ |
| $z_i^{p,d}$ | Binary variable if driver $d$ drops off participant $p$ at node $i$ |
| $\mu^p$ | Binary variable if participant $p$ is matched |
| $\delta^d$ | Binary variable if driver $d$ provide drive service |
| $\kappa_i^d$ | Number of people in vehicle $d$ after serving node $i$ |
| $u_i^d$ | Counter on node $i$ in the path of driver $d$ |
| $\tau_i^p$ | The time at which participant $p$ visits node $i$ |

In the notation provided, we formulate the many-to-many flex-role ridesharing problem within a network that includes HOV as follows.

$$\max \sum_{p} \lambda_p \mu^p - \sum_{i,j} \sum_{p,d} t_{i,j}\, x_{i,j}^{p,d} \tag{1}$$

$$\sum_{j \in \mathbb{AN}_d} x_{o_d,j}^{d,d} = \delta^d \qquad \forall\, d \in \mathbb{C} \tag{2}$$

$$\sum_{i \in \mathbb{AN}_d} x_{i,e_d}^{d,d} = \delta^d \qquad \forall\, d \in \mathbb{C} \tag{3}$$

$$\sum_{i \in \mathbb{AN}_d} x_{i,j}^{d,d} = \sum_{k \in \mathbb{AN}_d} x_{j,k}^{d,d} \qquad \forall\, d \in \mathbb{C}, \forall\, j \in \mathbb{AN}_d \setminus \{o_d, e_d\} \tag{4}$$

$$\sum_{d \neq p, d \in \mathbb{C}} \sum_{j \in \mathbb{AN}_{pd}} x_{o_p,j}^{p,d} = \mu^p \qquad \forall p \in \mathbb{P} \tag{5}$$

$$\sum_{d \neq p, d \in \mathbb{C}} \sum_{i \in \mathbb{AN}_{pd}} x_{i,e_d}^{p,d} = \mu^p \qquad \forall p \in \mathbb{P} \tag{6}$$

$$\sum_{d \neq p, d \in \mathbb{C}} \sum_{i \in \mathbb{AN}_{pd}} x_{i,j}^{p,d} = \sum_{d \neq p, d \in \mathbb{C}} \sum_{k \in \mathbb{AN}_{pd}} x_{j,k}^{p,d} \qquad \forall p \in \mathbb{P}, \forall\, j \in \mathbb{AN}_{pd} \setminus \{o_p, e_p\} \tag{7}$$

$$\sum_{d \neq p, p \in \mathbb{P}} x_{i,j}^{p,d} \leq Q_d\, x_{i,j}^{d,d} \qquad \forall\, d \in \mathbb{C}, (i,j) \in \mathbb{AE}_{pd} \tag{8}$$

$$\sum_{d \in \mathbb{C}} x_{i,j}^{p,d} \leq 1 \qquad \forall p \in \mathbb{P}, (i,j) \in \mathbb{AE}_{pd} \tag{9}$$

$$x_{i,j}^{p,d} \leq \delta^d \qquad \forall p \in \mathbb{P}, d \in \mathbb{C}, (i,j) \in \mathbb{AE}_{pd} \tag{10}$$

$$\delta^d + \mu^d \leq 1 \qquad \forall\, d \in \mathbb{C} \tag{11}$$

$$y_i^{p,d} \leq \sum_{j \in \mathbb{AN}_{pd}} x_{i,j}^{p,d} \qquad \forall p \in \mathbb{P}, d \in \mathbb{C}, i \in \mathbb{AN}_{pd} \quad (12)$$

$$2y_i^{p,d} - 1 \leq \sum_{j \in \mathbb{AN}_{pd}} x_{i,j}^{p,d} - \sum_{j \in \mathbb{AN}_{pd}} x_{j,i}^{p,d} \leq y_i^{p,d} \qquad \forall p \in \mathbb{P}, d \in \mathbb{C}, i \in \mathbb{AN}_{pd} \quad (13)$$

$$2z_i^{p,d} - 1 \leq \sum_{j \in \mathbb{AN}_{pd}} x_{j,i}^{p,d} - \sum_{j \in \mathbb{AN}_{pd}} x_{i,j}^{p,d} \leq z_i^{p,d} \qquad \forall p \in \mathbb{P}, d \in \mathbb{C}, i \in \mathbb{AN}_{pd} \quad (14)$$

$$\sum_{i \in \mathbb{AN}_{pd}} y_i^{p,d} = \sum_{j \in \mathbb{AN}_{pd}} z_j^{p,d} \qquad \forall p \in \mathbb{P}, d \in \mathbb{C} \quad (15)$$

$$u_j^d \geq u_i^d + x_{i,j}^{d,d} - |\mathbb{AE}_d|\left(1 - x_{i,j}^{d,d}\right) \qquad \forall d \in \mathbb{C}, (i,j) \in \mathbb{AE}_d \quad (16)$$

$$u_j^d \geq u_i^d + 1 - |\mathbb{AN}_d|\left(2 - z_j^{p,d} - y_i^{p,d}\right) \qquad \forall p \in \mathbb{P}, d \in \mathbb{C}, i,j \in \mathbb{AN}_{pd} \quad (17)$$

$$\kappa_j^d \leq \kappa_i^d + \sum_{p \in \mathbb{P}} y_j^{p,d} - \sum_{p \in \mathbb{P}} z_j^{p,d} + M_1\left(1 - x_{i,j}^{d,d}\right) \qquad \forall d \in \mathbb{C}, (i,j) \in \mathbb{AE}_d \quad (18)$$

$$\kappa_j^d \geq \kappa_i^d + \sum_{p \in \mathbb{P}} y_j^{p,d} - \sum_{p \in \mathbb{P}} z_j^{p,d} - M_1\left(1 - x_{i,j}^{d,d}\right) \qquad \forall d \in \mathbb{C}, (i,j) \in \mathbb{AE}_d \quad (19)$$

$$\kappa_i^d \geq n_H - Q^d\left(1 - x_{i,j}^{d,d}\right) \qquad \forall d \in \mathbb{C}, (i,j) \in \mathbb{AE}_d \cap \mathbb{E}_H \quad (20)$$

$$\tau_i^d + t_{i,j} \leq \tau_j^d + M_{i,j}^d\left(1 - x_{i,j}^{d,d}\right) \qquad \forall d \in \mathbb{C}, (i,j) \in \mathbb{AE}_{pd} \quad (21)$$

$$ED_d(\delta^d + \mu^d) \leq \tau_i^d \leq LA_d(\delta^d + \mu^d) \qquad \forall d \in \mathbb{C}, \forall i \in \mathbb{AN}_d \quad (22)$$

$$ED_p\, \mu^p \leq \tau_i^p \leq LA_p \mu^p \qquad \forall p \in \mathbb{P}, \forall i \in \mathbb{AN}_p \quad (23)$$

$$x_{i,j}^{p,d}, y_j^{p,d}, z_i^{p,d}, \mu^p, \delta^d \in \{0,1\} \qquad \forall d \in \mathbb{C}, \forall p \in \mathbb{P}, i,j \in \mathbb{AN}_{pd} \quad (24)$$

$$\kappa_i^d, u_i^d \in \mathcal{N}, \tau_i^p \in Z^+ \qquad \forall d \in \mathbb{C}, \forall p \in \mathbb{P} \quad (25)$$

Objective function (1) maximizes total matches and minimizes travel time, with $\lambda_p$ as the travel time for the fastest route from origin to destination, and $t_{i,j}$ as edge travel time. Constraint (2) ensures the departure of drivers from their origins, while constraint (3) guarantees their arrival at their destinations. Constraint (4) fulfills flow conservation for the driver's path. Constraint (5) directs riders out of their origin and (6) ensures that they end their trips at their destination. Constraint (7) indicates flow conservation for riders. Moreover, these constraints enable passengers to transfer between drivers, creating a many-to-many scheme. Constraint (8) ensures vehicle capacity and the prerequisite that a driver must travel along a link for a rider on that link. Constraint (9) guarantees single-driver assignments for each step of a rider's path. Constraint (10) mandates car owners to be assigned as drivers before providing service to others. Constraint (11) states that car owners can either provide a ride service or receive it. Constraint set (12) states that a pick follows and immediate travel to the next node. Constraint (13) checks the feasibility of pickup at node $i$ by driver $d$. This requires that incoming edges at $i$ sum to zero, while outgoing edges sum to one. Constraint (14) monitors the possibility of drop-offs, ensuring that a source node only has outgoing flow while a destination node only receives incoming flow. Constraints (15)-(17) ensure meaningful orders for pickups and drop-offs. Constraint (15) guarantees that each pickup corresponds to a drop-off. Constraint (16) act as the counter on nodes visited by driver $d$ as they cross edges, tracking progress along the route. For pickup-drop-off order, constraint (17) ensures pick up location $i$ have a higher counter than drop-off point $j$. Constraint sets (18) and (19) limit the number of passengers on board after serving node $j$. The coefficient $M_1$ is set to $M_1 = Q^d + \sum_{p \in \mathbb{P}} 1$. Constraints (20) enforces a minimum passenger requirement for HOV lane access. Constraint (21) relates departure times at nodes and its subsequent immediate node, where $M_{i,j}^p = LA_p - ED_p + t_{i,j}$. This also prevents subtours. Constraint (22) and (23) ensures adherence to time windows for participants. Constraint sets (24) and (25) describe the type of variables.

# 3 Methodology

To compare the time savings from HOV and Flex roles, we utilize the New York City Taxi Trip dataset [7], which contains coordinates, timestamps, and trip durations. We randomly selected ~17k requests, with ~8k from non-car owners and ~9k from car owners. To address missing departure and arrival time in data, we follow literature [1] practices by adding or subtracting a random time between 0 and 5 minutes from the pickup and drop-off times. To identify accessible nodes within a user's time frame, we employ Yen's shortest path algorithm [8] and analyze the first 10 shortest paths. We obtain map data using the *OpenStreetMapX* package [9]. The map contains 40490 nodes and 12,242 edges, from which we randomly set 785 edges to have HOV lanes (see Fig.1). A minimum of 2 occupants is required for HOV lane access. The system operates from 9:00 am to 07:00 pm, and we employ a rolling horizon [5] with a 5-minute period for all experiments. All models are implemented in Julia v1.8.3 and solved using ILOG CPLEX v22.1.0 on a laptop with a Core (TM) i7 processor and 16GB of RAM.

Fig.1 HOV lanes in the Manhattan area [2]

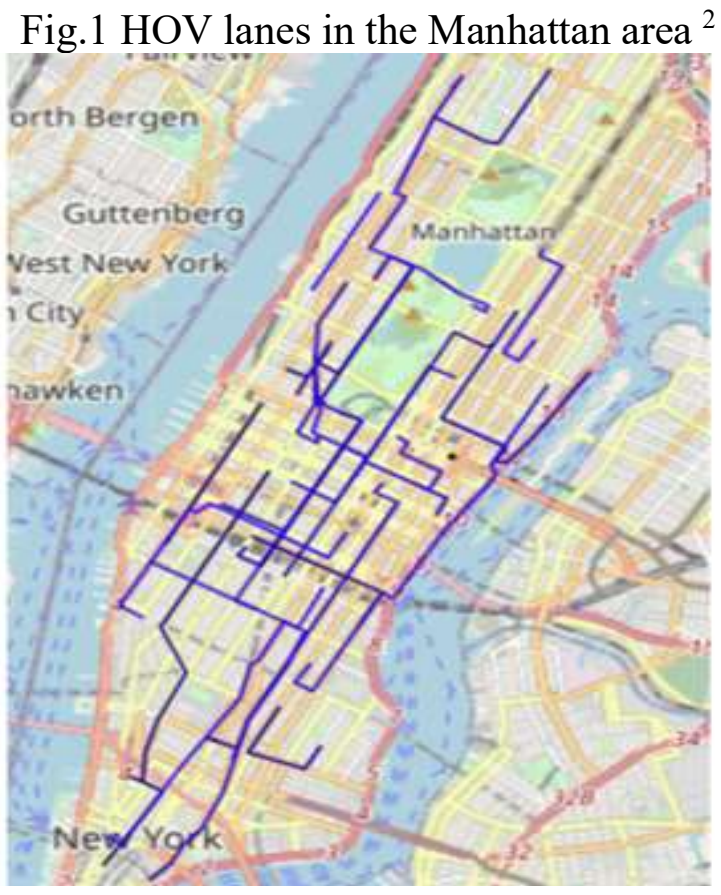

# 4 Results

Table 2 summarizes findings, with values averaging over the horizon. The '%Match' represents served participants to non-car owners, while '%Shift' is the ratio of car owners served by others. Combining role flexibility with HOV lanes achieves the lowest objective function. However, role flexibility alone achieves a higher %match, likely due to HOV access requirements that limits combination of users can travel together. Role flexibility enhances total matches compared to inflexible models, but it is influenced by the driver-passenger ratio. We sometimes achieve 112% matches, serving all non-car owners and 12% of car owners through role shifts. However, in other cases, more role shifts can lead to imbalances, affecting service levels to non-car owners.

[2] Image produced with OpenStreetMapX [9]; MIT licensed.

# 5 Conclusion

Role flexibility improves total matches in both HOV and non-HOV scenarios, with its effectiveness tied to the driver-passenger ratio. Further research on the correlations among role flexibility, role shifts, participant ratios, and service levels is an interesting avenue for future research.

**Table 2.** Effect of HOV and Role Flexibility on Matching

| ***Model*** | ***Objective*** | %***Match*** | %***Role Shift*** |
|---|---|---|---|
| Flex Role-HOV | 1164.13 | 80.62 | 27.87 |
| Flex Role-No HOV | 1489.12 | 89.12 | 29.80 |
| No Flex Role-HOV | 2666.74 | 60.14 | - |
| No Flex Role-No HOV | 3153.18 | 59.41 | - |